\documentclass[12pt]{article}
\usepackage{amsmath,amssymb,amstext}
\usepackage{amsthm}

\newcommand{\Li}{\mbox{Li}}

\newtheorem{definition}{Definition}
\newtheorem{theorem}{Theorem}

\newtheorem{cor}{Corollary}
\newtheorem{lemma}{Lemma}[section]

\newcommand{\C}{\mathbb{C}}
\newcommand{\Q}{\mathbb{Q}}
\newcommand{\R}{\mathbb{R}}
\newcommand{\Z}{\mathbb{Z}}
\newcommand{\p}{\mathfrak {p}}
\newcommand{\Pri}{\mathfrak {P}}

\newcommand{\Gal}{\operatorname{Gal}}

\begin{document}
\title{\bf Chebotarev Sets}
\author{Hershy Kisilevsky,\ Michael O. Rubinstein
\footnote{Support for work on this paper was provided by the
National Science Foundation under awards DMS-0757627 (FRG grant),
and NSERC Discovery Grants}}

\maketitle

\abstract{We consider the problem of determining whether a set of primes, or,
more generally, prime ideals in a number field, can be realized as a finite
union of residue classes, or of Frobenius conjugacy classes. We give
necessary conditions for a set to be realized in this manner, and show that the subset of
primes consisting of every other prime cannot be expressed in this way, even
if we allow a finite number of exceptions.}

\baselineskip=18pt
\parskip=13pt

\setcounter{page}{1}

\section{Introduction}

In this paper, we consider the following problem, and its generalizations:

Let $p_n$ denote the $n$-th prime, $\pi$ the set of all primes, and
$P_{\text{odd}}$ the set consisting of every other prime:
\begin{equation}
    \label{eq:P}
    P_{\text{odd}} = \{p_n \in \pi | n \text{ odd}\}
    =\{2,5,11,17,23,\ldots\}.
\end{equation}

Can the set $P_{\text{odd}}$ be realized as a finite union of primes in residue
classes, even if we are willing to allow a finite number of exceptions?
 
More generally, can we realize a given set of primes ideals in a number field $K$ with rational 
density relative to the full set of primes ideals of $K$ as a finite union of prime ideals that arise in the
Chebotarev density theorem, i.e. of Frobenius conjugacy classes in the
Galois groups of finite Galois extensions of $K$?

In the following we will identify the (positive) prime number $p$ with the ideal that it generates in the 
ring of integers.

The natural instinct is that the set $P_{\text{odd}}$ above cannot be realized
in this manner. In fact, the cardinality of subsets of primes with density $1/2$ is
uncountable, for example we can pick subsets by performing a random coin flip at
each prime, yet the number of residue classes (or, more generally, Frobenius
conjugacy classes) is countable. Therefore, most subsets of primes will
fail to arise from residue or conjugacy classes.

In this article we prove that $P_{\text{odd}}$  cannot
be realized as a finite union of primes in residue/conjugacy classes,
even if we are willing to allow a finite number of exceptions.

We will show that the set $P_{\text{odd}}$ is too `quiet' relative to the
set of all primes to arise from such arithmetic sets. Primes in
progressions, or in Frobenius conjugacy classes are quantifiably irregular,
as a result of the non-trivial zeros of the $L$-functions that govern them.
 
\subsection{Chebotarev Sets}
\label{sec:Chebotarev intro}

Let $K$ be a number field, and let $\pi(K)$ denote the set of all non-zero prime
ideals of $K$. Let $L/K$ be a finite Galois extension of number fields with
Galois group $\Gal(L/K)=G$. For a prime ideal $\p \in \pi(K),$ unramified in $L/K$,
let $(L/K,\p)$
denote the conjugacy class of Frobenius automorphisms
of primes $\Pri \in \pi(L),$ with $\Pri$ dividing $\p.$
If $C \subset G$ is a conjugacy class in $G$, let
$\pi(L/K,C)$ denote the set of primes $\p \in \pi(K)$ such that $\p$ is unramified
in $L/K$ and such that $(L/K,\p)$ is equal to $C$.

For two sets $S_1,S_2$ the symmetric difference $S_1 \bigtriangleup S_2$ is the
set $S_1 \setminus S_2 \cup S_2 \setminus S_1$. We will say two sets $S_1$ and $S_2$ are equal 
up to finite sets if $S_1 \bigtriangleup S_2$ is finite.

\begin{definition}
    Call a set
    of primes $P \subseteq \pi(K)$ a {\it Chebotarev set\/} for $K$ if there are
    finitely many finite Galois extensions $L_i/K$ and conjugacy classes $C_i
    \subset \Gal(L_i/K)$ such that $P= \cup_i \pi(L_i/K,C_i)$ up to finite sets.
    That is  $P \bigtriangleup\cup_i \pi(L_i/K,C_i)$ is finite. 
\end{definition}

We allow the possibility that $P$ is the empty set (i.e. that the number of fields $L_i$ is zero).

Suppose that $K \subseteq L' \subseteq L$ is a tower of finite extensions with
$L'/K$ and $L/K$ both Galois. Then a conjugacy class in $\Gal(L'/K)$ can be
lifted to a finite union of conjugacy classes in $\Gal(L/K)$. Hence $P$ is a
Chebotarev set for $K$ if and only if there is a finite Galois extension $L/K$
and a finite number of distinct conjugacy classes $C_i \subset \Gal(L/K)$ such that $P=
\cup_i \pi(L/K,C_i)$ up to finite sets. 

Given Chebotarev sets $P,Q \subseteq \pi(K)$ the following properties are consequences of  the Chebotarev
density theorem:\\
$\bullet$  the natural and Dirichlet densities $\delta(P)$ exist and are equal to
$$
\frac{\sum_i \vert C_i \vert }{ \vert G \vert }
$$   \\
$\bullet$  if $\delta(P)=0$ or $1$ then $P$ is finite or co-finite \\
$\bullet$   if $P,Q$ are Chebotarev sets then so are $P \cup Q, P \cap Q, \pi(K) \setminus P$,
$P \Delta Q$.

 Thus, any
subset of $\pi(K)$ with no density or irrational density cannot be a Chebotarev
set. Also, any set of density $0$ (resp. $1$) which is infinite (resp. co-infinite) cannot be a Chebotarev set.
So, for example, the set of (positive) primes    congruent to $3 \pmod 4$ constitute a Chebotarev set of $\Q$
but they are generators of an infinite set of prime ideals of density $0$ when viewed in $K=\Q(i)$ and 
therefore do {\em not} form a Chebotarev set of $K.$ 

Our problem is to produce a set $P \subset \pi(\Q)$ of positive rational density which is
provably not a Chebotarev set. 

For a subset $P \subset \pi(K)$ of prime ideals of $K$, let $P(x)$ denote the function which 
counts the number of elements in $P$ and with absolute norm $\leq x$:
\begin{equation}
    \label{eq:P(x) C}
    P(x)=\#\{\p \in P | N_{K/\Q}(\p) \leq x\}.
\end{equation}
We also let
\begin{equation}
    \pi(x,K)=\#\{\p \in \pi(K) | N_{K/\Q}(\p) \leq x\}.
\end{equation}

To this end, we have the following theorem.

\begin{theorem}
\label{thm:1}
Let $\beta \in \Q$, with $0 < \beta < 1$.
Assume that $P \subset \pi(K)$ is a Chebotarev set of density $\beta.$
Then,
\begin{equation}
    \label{eq:P(x) K vs pi(x) K}
    P(x) - \beta \pi(x,K) =
    \Omega\left(\frac{x^{1/2}}{\log{x}}\right).
\end{equation}
with the implied constant in the $\Omega$ depending on $P$.
\end{theorem}

Here, for real functions $f,g>0$, with $g>0$,
we are using the $\Omega$ notation to mean the following:
We say that $$f(x) = \Omega(g(x))$$ if $$\lim\sup|f(x)|/g(x) >0,$$ i.e. if
there is a constant $c>0$ and a sequence ${x_n}$ with $x_n\to\infty$, such that
$|f(x_n)| > c g(x_n)$.

\begin{cor}
\label{cor:1}
The set $P_{\text{odd}}$ is not a Chebotarev set.
\end{cor}
\begin{proof}
The set $P_{\text{odd}}$ specified in the introduction has natural density (and hence Dirichlet density)  $1/2$ in $\pi(\Q).$
Because $P_{\text{odd}}$ consists of every other prime, we have that
\begin{equation}
    \label{eq:P(x) vs pi(x)}
    P_{\text{odd}}(x) -  \pi(x)/2 =
    \begin{cases}
        1/2,\quad \text{if $p_{2j-1}\leq x< p_{2j}$} \\
        0,\quad \text{if $p_{2j}\leq x< p_{2j+1}$}.
    \end{cases}
\end{equation}
Therefore Theorem 1 implies that $P_{\text{odd}}$ cannot be realized as a Chebotarev set.
\end{proof}

As we will explain in the proof of this theorem, the same conclusion  holds if
we replace the counting function $\beta \pi(x,K)$ by
$\beta \Li(x)$, or with the counting function $sQ(x)$
of another Chebotarev set $Q \subseteq \pi(K)$ of non-zero density $\delta(Q)$, which is
essentially distinct from $P$, so long as the density, $s \delta(Q)=\beta$. By
essentially distinct, we mean that the symmetric difference $P \bigtriangleup
Q$ is infinite.

\begin{theorem}
\label{thm:1b}
Let $P \subset \pi(K)$ be a Chebotarev  set of density $\beta,$ $0 <\beta <1,$ and let
$f(x)$ stand for $\beta \Li(x)$, or $s Q(x)$, as  above.
Then
\begin{equation}
    \label{eq:P(x) K vs pi(x) F}
    P(x) - f(x) =
    \Omega\left(\frac{x^{1/2}}{\log{x}}\right),
\end{equation}
with the implied constant depending on $P$ and $f$.
\end{theorem}

In a personal communication with the authors, Serre raised the question of
whether our techniques allow us to address the size of
a summatory function on prime ideals of $K$ with complex valued weight function that is
constant on Frobenius conjugacy classes of $G$, i.e. a sum of the form
\begin{equation}
    \label{eq:serre}
    \sum_{\p: N\mathfrak{p}\leq X \atop \p \text{unramified}} h(\p)
\end{equation}
with $h$ complex valued and taking on finitely many values according to
conjugacy classes of $G$. Such a summatory function can be expressed as a linear
combination of counting functions $\pi(x,L/K,C)$.

More specifically,
let $C_1, \ldots, C_r$ be the distinct
conjugacy classes of $G$. Let $\eta_1, \ldots, \eta_r$ be complex numbers not all
$0$, let $\eta$ be defined by
\begin{equation}
    \label{eq:lambda_0}
    \eta := \sum_{j=1}^r \eta_j \delta(C_j)= \frac{1}{\vert G
    \vert}\sum_{j=1}^r \eta_j \vert C_j \vert,
\end{equation}
and let $\pi(x,L/K,C)$ be defined by~\eqref{eq:pi C}.

\begin{theorem}
\label{thm:1c}
Let
\begin{equation}
    \label{eq:F}
    F(x) = \sum_{j=1}^r \eta_j \pi(x,L/K,C_j).
\end{equation}
Then
\begin{equation}
    \label{eq:F vs Li}
    F(x) - \eta \Li(x) =
    \Omega\left(\frac{x^{1/2}}{\log{x}}\right),
\end{equation}
with the implied constant depending on the choice of the set $\{\eta_j\}$. Furthermore, the same result holds if we replace
$\Li(x)$ by $\pi(x,K)$, though one then needs the additional restriction that not all the $\eta_j$ are equal.
\end{theorem}

The key idea used to prove these theorems is that the functions on the left hand side 
of~\eqref{eq:P(x) K vs pi(x) K} \eqref{eq:P(x) K vs pi(x) F}, and
\eqref{eq:F vs Li} are discontinuous at an infinite number of values of $x.$ 
Consequently, when expressed as a linear combination of
explicit formulas, infinitely many of the non-trivial zeros of the relevant
$L$-functions must survive.
These zeros (in fact we only need one non-trivial zero to enter) are responsible for
making these differences large on average,
which we show by considering their mean square on a logarithmic scale.

Note that the statements of Theorems~\ref{thm:1}-~\ref{thm:1c} do not assume
the Generalized Riemann Hypothesis. In fact, if the GRH does not hold,
stronger $\Omega$ results than these hold, hence we have
stated these theorems unconditional on the GRH.

The precise statement of the $\Omega$ bound in the case that the GRH fails
requires some discussion concerning the location of the zeros of the relevant
$L$-functions, and how these zeros interact upon taking certain linear
combinations of the logarithmic derivatives of these $L$-functions. This
discussion and corresponding result can be found in Section~\ref{sec:theta >
1/2} and in Theorem~\ref{thm:3} at the end of Section~\ref{sec:proofs in
general}.

Also observe, in our theorems we do not prove $\Omega_{\pm}$ results, i.e. we do
not address the question of sign changes. Without further assumptions, such as linear
independence of the non-trivial zeros over $\Q$ (other than those possibly occurring at $s=1/2$),
we cannot, in general, prove the existence of sign changes.
For a discussion on issues related to sign changes see~\cite{RS}.

\begin{definition}
    Call a set of primes $P \subseteq \pi(K)$ an {\it almost Chebotarev set\/}
    for $K$ if there is a Chebotarev set $Q \subseteq \pi(K)$ such that $P = Q$ up to sets
    of density zero. That is  $P \bigtriangleup Q$ has density zero.
\end{definition}

It seems much more difficult to prove the existence of a set which is  not almost
Chebotarev -- although we suspect that our example, $P_{\text{odd}}$, is one such set.

We conclude the introduction by noting that Serre studied `Frobenian' (i.e.
named differently than here) sets and functions in his paper~\cite{S} and
book~\cite{S2} (see his Chapter 3), the latter in relation to the problem of counting the number of
solutions mod $p$ to a system of polynomial equations. Lagarias defined a
similar notion of Chebotarev sets in~\cite{L}, also for studying solutions to
polynomial congruences modulo $p$. See Lemma 3.1 in his paper for the
equivalence of his definition to ours, though without allowing for finitely
many exceptions.

\section{The classical case}
\label{sec:classical}

In this section we consider the more classical situation of sets $P$ of
rational primes that are realized using residue classes. We will essentially
establish Theorem~\ref{thm:1} for the special case of residue classes, rather
than Frobenius conjugacy classes.

The techniques that we develop will serve as a model, in Section~\ref{sec:chebotarev},
where we will modify our approach to the general setting of Chebotarev sets.

Assume that
\begin{equation}
    \label{eq:P as union}
    P = P_0 \cup \bigcup_{j=1}^r \pi(q_j,a_j) \setminus P_1
\end{equation}
where $P_0,P_1$ consists of finitely many elements, i.e. the possible exceptions in excess and deficiency,
and
\begin{equation}
    \label{eq:P a,q}
    \pi(q,a) = \{p \text{ prime} | p = a \mod q\}
\end{equation}
consists of rational primes in the residue class $a$ mod $q$.


 As noted in the introduction there is a single positive integer $q$,
and distinct residue classes $a_j$ mod $q$ such that, after relabelling $r$
as needed,
\begin{equation}
    \label{eq:P as union b}
    P = P_0 \cup \bigcup_{j=1}^r \pi(q,a_j) \setminus P_1.
\end{equation}
We can also assume that $\gcd(a_j,q)=1$,
that $P_0$ is disjoint from
$\bigcup_{j=1}^r \pi(q_j,a_j)$, and $P_1$ is contained within this union.
 
Next, we define, as usual,
\begin{equation}
    \pi(x,q,a) = \# \{p \leq x \mid p = a \bmod q\}.
\end{equation}
From equation~\eqref{eq:P as union b}, we
have, for $x$ larger than all of the elements of $P_0$ and $P_1$, that
\begin{equation}
    \label{eq:P from pi}
    P(x) = \lambda + \sum_{j=1}^r \pi(x,q,a_j),
\end{equation}
where
\begin{equation}
    \label{eq:lambda P}
    \lambda= |P_0| - |P_1|.
\end{equation}

The prime number theorem states that
\begin{equation}
    \pi(x) \sim \Li(x),
\end{equation}
where
\begin{equation}
    \Li(x) = \int_2^x dt/\log{t}
    \sim \frac{x}{\log{x}}.
\end{equation}
The prime number theorem for arithmetic progressions, proven
by Hadamard and de la Vall\'{e} Poussin,
asserts, for $\gcd(a,q)=1$, that primes are equidistributed amongst
the residue classes mod $q$ that are relatively prime to $q$:
\begin{equation}
    \label{eq:pnt progressions}
    \pi(x,q,a) \sim \frac{\pi(x)}{\phi(q)},
\end{equation}
where $\phi(q)= \#\{ a \bmod q \mid (a,q)=1 \}$.
Therefore,
\begin{equation}
    \label{eq:P vs pi}
    P(x) \sim \frac{r}{\phi(q)} \pi(x).
\end{equation}
We also make the assumption that
\begin{equation}
    \frac{r}{\phi(q)} < 1,
\end{equation}
so that $P$ omits a positive proportion of the primes.

\subsection{Outline of proof}

We now provide a summary of our proof of Theorem~\ref{thm:1} in the case of
residue classes. We begin with the classical formula~\eqref{eq:psi progressions}
which uses Dirichlet characters to count prime powers in arithmetic progressions,
weighted by the von-Mangoldt function. A summation by parts, and
isolating the squares and higher powers of primes, allows us to pass to
counting primes in arithmetic progressions in formula~\eqref{eq:pi progressions
in terms of Pi c}. This leads to formula~\eqref{eq:main formula} for $P(x)
-r\pi(x)/\phi(q)$ in Lemma~\ref{lemma:P(x) formula}.

We then apply the explicit formula, in Section~\ref{sec:explicit formula}, to
derive equation~\eqref{eq:main formula zeros} in Lemma~\ref{lemma:P explicit}.
That Lemma expresses $P(x)-r\pi(x)/\phi(q)$ in terms of non-trivial zeros of
Dirichlet $L$-functions, along with two complicated, but innocuous, terms which
are denoted by $A(x)$ and $B(x)$. We then exploit the fact that
$P(x)-r\pi(x)/\phi(q)$ is discontinuous at the primes to show, in
Lemma~\ref{lemma:alpha_rho nonzero}, that the sum over zeros is non-empty.

Having at least one non-trivial zero appear in~\eqref{eq:main formula zeros}
i.e. with $\alpha_\rho\neq 0$,
guarantees the $\Omega$ bound of Theorem~\ref{thm:1}, whether the
Generalized Riemann Hypothesis is false or not. The former case is
considered in Theorem~\ref{thm:theta > 1/2} of Section~\ref{sec:theta > 1/2}.
Two proofs in the latter case is handled in Sections~\ref{sec:theta = 1/2}-
\ref{sec:mean square b} by considering mean squares of $P(x)-r\pi(x)/\phi(q)$.

Finally, we generalize our results from the classical case of residue classes to
Chebotarev sets in Section~\ref{sec:chebotarev}.

\subsection{Counting primes in arithmetic progressions}
\label{sec:arithmetic progressions}

Our first step is to derive the following Lemma, which expresses
$P(x)-r\pi(x)/\phi(q)$ in terms of Dirichlet characters.

\begin{lemma}
\label{lemma:P(x) formula}
For $x$ larger than all the elements of $P_0$ and $P_1$, we have
\begin{eqnarray}
    \label{eq:main formula}
    &&P(x) - \frac{r}{\phi(q)}\pi(x) =
    \lambda + \sum_{j=1}^{r} \pi(x,q,a_j) -\frac{r}{\phi(q)}\pi(x) \notag \\
    &=& \lambda +
    \sum_{\chi \bmod q \atop \chi\neq\chi_0} c_\chi
    \left(
        \frac{\psi(x,\chi)}{\log(x)}
        +\int_2^x \frac{\psi(t,\chi)}{t \log(t)^2} dt
    \right)
    +\frac{r R(x,1)}{\phi(q)}
    -\sum_{j=1}^{r}
    R(x,q,a_j) \notag \\
    &-&\frac{r}{\phi(q)}
    \left(
         \frac{\sum_{p^k \leq x \atop  p \mid q} \log(p)}{\log(x)}
        +\int_2^x \frac{ \sum_{p^k \leq t \atop  p \mid q} \log(p) }{t \log(t)^2} dt
    \right),
\end{eqnarray}
where $\lambda$ is defined in~\eqref{eq:lambda P},
$R(x,1)$ and $ R(x,q,a)$ are defined, below, in~\eqref{eq:R(x,1)}
and~\eqref{eq:remainder},
and
\begin{equation}
    \label{eq:c_chi}
    c_\chi = \frac{1}{\phi(q)} \sum_{j=1}^{r} \bar{\chi}(a_j).
\end{equation}
\end{lemma}

\begin{proof}

For a given $q$, let $\chi$ be a Dirichlet character modulo $q$.
We will denote the principal character by $\chi_{0}$.

As in the proof of Dirichlet's Theorem, it is easier to count prime powers
weighted by the von-Mangoldt function than it is to simply count primes. Thus, let
\begin{equation}
    \psi(x,\chi) := \sum_{n \leq x} \chi(n) \Lambda(n),
\end{equation}
where $\Lambda(n)= \log p$ if $n=p^{m}$ for some $m \in \Z$, and is 0
otherwise.

We have
\begin{eqnarray}
    \psi(x,q,a) & := & \sum_{{n \leq x}\atop{n \equiv a \bmod q}} \Lambda(n) \notag \\
    & = & \frac{1}{\phi(q)}\sum_{\chi \bmod q} \bar{\chi}(a)
    \sum_{n \leq x}\Lambda(n)\chi(n) \notag \\
     & = &
    \frac{1}{\phi(q)} \sum_{\chi \bmod q} \bar{\chi}(a) \psi(x,\chi).
    \label{eq:psi progressions}
\end{eqnarray}
The main contribution to $\psi(x,q,a)$ comes from the principal character:
\begin{equation}
    \label{eq:principal psi}
    \psi(x,\chi_0) = \sum_{p^k \leq x \atop  p \nmid q} \log(p) = \psi(x) - \sum_{p^k \leq x \atop  p | q} \log(p).
\end{equation}
Therefore
\begin{equation}
    \psi(x,q,a)
    = \frac{1}{\phi(q)}
    \left(
        \psi(x)
        + \sum_{\chi \bmod q \atop \chi\neq\chi_0} \bar{\chi}(a) \psi(x,\chi)
        - \sum_{p^k \leq x \atop  p \mid q} \log(p)
    \right).
    \label{eq:psi progressions b}
\end{equation}

We define
\begin{eqnarray}
    \Pi(x,q,a)
    & := & \sum_{{n \leq x}\atop{n \equiv a \bmod q}} \frac{\Lambda(n)}{\log(n)}
    =  \sum_{{p^k \leq x}\atop{p^k \equiv a \bmod q}} \frac{1}{k} \notag \\
    & = & \pi(x,q,a) + R(x,q,a)
\end{eqnarray}
with
\begin{equation}
    \label{eq:remainder}
    R(x,q,a) =
    \sum_{{p^k \leq x\atop k\geq 2}\atop{p^k \equiv a \bmod q}} \frac{1}{k}.
\end{equation}
Therefore,
\begin{equation}
    \label{eq:pi progressions in terms of Pi}
    \pi(x,q,a) = \Pi(x,q,a) - R(x,q,a).
\end{equation}
Now, summing by parts
\begin{equation}
    \Pi(x,q,a) = \frac{\psi(x,q,a)}{\log(x)} + \int_2^x \frac{\psi(t,q,a)}{t \log(t)^2} dt,
\end{equation}
so that
\begin{equation}
    \label{eq:pi progressions in terms of Pi b}
    \pi(x,q,a) =
     \frac{\psi(x,q,a)}{\log(x)} + \int_2^x \frac{\psi(t,q,a)}{t \log(t)^2} dt
    -R(x,q,a).
\end{equation}
Substituting~\eqref{eq:psi progressions b} into the above,
we get
\begin{eqnarray}
    \label{eq:pi progressions in terms of Pi c}
    &&\pi(x,q,a) = \notag \\
    &&\frac{1}{\phi(q)}
    \left(
        \frac{\psi(x)}{\log(x)} +
        \int_2^x \frac{\psi(t)}{t \log(t)^2} dt
        - \frac{\sum_{p^k \leq x \atop  p \mid q} \log(p)}{\log(x)}
        -\int_2^x \frac{ \sum_{p^k \leq t \atop  p \mid q} \log(p) }{t \log(t)^2} dt
    \right) \notag \\
    &+&
    \frac{1}{\phi(q)}
    \sum_{\chi \bmod q \atop \chi\neq\chi_0} \bar{\chi}(a)
    \left(
        \frac{\psi(x,\chi)}{\log(x)}
        +\int_2^x \frac{\psi(t,\chi)}{t \log(t)^2} dt
    \right)
     -R(x,q,a). \notag \\
\end{eqnarray}
The special case of $q=1$ (and any value for $a$) of~\eqref{eq:pi progressions in terms of Pi b} is:
\begin{eqnarray}
    \label{eq:pi in terms of Pi}
    \pi(x) =
     \frac{\psi(x)}{\log(x)} + \int_2^x \frac{\psi(t)}{t \log(t)^2} dt
     -R(x,1),
\end{eqnarray}
with
\begin{equation}
    \label{eq:R(x,1)}
    R(x,1) =
    \sum_{k\geq 2} \frac{\pi(x^{1/k})}{k}.
\end{equation}

With these formulas in hand, consider again the difference $P(x) - \frac{r}{\phi(q)}\pi(x)$.
Subtracting $\frac{r}{\phi(q)}\pi(x)$ from~\eqref{eq:P from pi}, and
then substituting~\eqref{eq:pi progressions in terms of Pi c} and \eqref{eq:pi in terms of Pi}
gives the Lemma.

Note, that, on summing over $r$ values of $a_j$ and
subtracting $r/\phi(q)$ times~\eqref{eq:pi in terms of Pi}, we cancelled the main term
\begin{equation}
    \label{eq:terms that cancel}
    \frac{r}{\phi(q)}\left(
        \frac{\psi(x,\chi_0)}{\log(x)}
        +\int_2^x \frac{\psi(t,\chi_0)}{t \log(t)^2} dt
    \right).
\end{equation}

\end{proof}

\subsection{Applying the explicit formula}
\label{sec:explicit formula}

We will show that
the right hand side of~\eqref{eq:main formula} can get as large, in absolute value, as
$\gg x^{1/2}/\log(x)$. This will be independent of the GRH.
In fact, if GRH fails, the lower bound that we can prove is at least as large.

To accomplish this, we will write an explicit formula for~\eqref{eq:main formula}
in terms of the zeros of the Dirichlet $L$-functions, $L(s,\chi)$, where
$\chi$ runs over all Dirichlet characters for the modulus $q$.

The explicit formula for $\psi(x,\chi)$, where $\chi\neq \chi_0$
is a primitive character, takes the form, for $x>1$ not a prime power,
\begin{equation}
    \label{eq:explicit formula}
    \psi(x,\chi) = - \sum_{\rho_\chi} \frac{x^{\rho_\chi}}{\rho_\chi}
    -(1-\mathfrak{a}_\chi) \log{x} -b(\chi) + \sum_{m=1}^\infty \frac{x^{\mathfrak{a}_\chi-2m}}{2m-\mathfrak{a}_\chi}, 
\end{equation}
where $\rho_\chi$ runs over the non-trivial zeros of $L(s,\chi)$, with the sum over
zeros taken as $\lim_{X \to \infty} |\Im{\rho_\chi}| < X$. Also, $\mathfrak{a}_\chi=1$ if
$\chi(-1)=-1$ and $0$ otherwise, and $b(\chi)$ is a constant depending on
$\chi$, namely the constant term in the Laurent expansion about $s=0$ of
$L'/L(s,\chi)$ (Taylor expansion if $\chi(-1)=-1$). If $x$ is a prime power,
the right hand side above converges to $\psi(x,\chi) - \Lambda(x) \chi(x)/2$, i.e. one
needs to subtract half of the last term in the sum defining $\psi(x,\chi)$. For
a derivation of this explicit formula see, for instance, Davenport \cite[pgs
115-120]{Davenport}.

When $\chi$ is an imprimitive character, say induced by $\chi_1$ mod $q_1$,
then
\begin{equation}
    \label{eq:psi imprimitive}
    \psi(x,\chi) = \psi(x,\chi_1) - \sum_{p^k \leq x \atop  p | q} \log(p) \chi_1(p^k).
\end{equation}
For notational convenience, in the case of imprimitive $\chi$, we set
$\mathfrak{a}_\chi = \mathfrak{a}_{\chi_1}$ and also $b(\chi) = b(\chi_1)$.

Therefore, we have shown that we can rewrite~\eqref{eq:main formula} in the following form.

\begin{lemma}
\label{lemma:P explicit}
For $x>\max(\lambda,1)$, and not a prime power,
\begin{eqnarray}
    \label{eq:main formula zeros}
    &&P(x) - \frac{r}{\phi(q)}\pi(x)
    = \frac{1}{\log{x}} \sum_{\rho} \alpha_\rho \frac{x^{\rho}}{\rho}
    + A(x)
    + B(x),
\end{eqnarray}
where $\alpha_\rho \in \C$ (described below in~\eqref{eq:alpha_rho}),
and the sum over $\rho$ is taken over the union over
the non-trivial zeros of all $L(s,\chi)$, $\chi$ mod $q$, $\chi\neq\chi_0$.

Here,
the function $A(x)$ gathers together all the
remaining terms that are discontinuous:
\begin{eqnarray}
    \label{eq:A(x)}
    A(x) &=& \frac{r R(x,1)}{\phi(q)}
    -\sum_{j=1}^{r}
    R(x,q,a_j)
    -\frac{r}{\phi(q)} \frac{\sum_{p^k \leq x \atop  p \mid q} \log(p)}{\log(x)} \notag \\
    &-&\frac{1}{\log{x}} \sum_{ {\chi \mod q\atop \chi \neq \chi_0} \atop \chi \text{ imprimitive}}
    c_\chi \sum_{p^k \leq x \atop  p | q} \log(p) \chi_1(p^k),
\end{eqnarray}
and $B(x)$ incorporates the rest:
\begin{eqnarray}
    \label{eq:B(x)}
    B(x) &=& \lambda
    -\frac{r}{\phi(q)} \int_2^x \frac{ \sum_{p^k \leq t \atop  p \mid q} \log(p) }{t \log(t)^2} dt
    +
    \sum_{\chi \bmod q \atop \chi\neq\chi_0}
    c_\chi \Big(\int_2^x \frac{\psi(t,\chi)}{t \log(t)^2} dt \notag \\
    &-&
    (1-\mathfrak{a}_\chi) + \frac{b(\chi)}{\log{x}}
    - \frac{1}{\log{x}}\sum_{m=1}^\infty \frac{x^{\mathfrak{a}_\chi-2m}}{2m-\mathfrak{a}_\chi}
    \Big).
\end{eqnarray}

\end{lemma}
Recall that $c_\chi$ is given in \eqref{eq:c_chi}. The role of the sum over imprimitive characters
in $A(x)$ is to account for~\eqref{eq:psi imprimitive}.

We assume that the $\rho$ in~\eqref{eq:main formula zeros} are distinct, by
grouping equal $\rho$ under the same $\alpha_\rho$. In the case of
imprimitive characters, the non-trivial zeros of $L(s,\chi)$ coincide with
those of the Dirichlet $L$-function, $L(s,\chi_1)$, corresponding to the
inducing character $\chi_1$. More precisely,
for given $\rho$,
\begin{equation}
    \label{eq:alpha_rho}
    \alpha_{\rho} = - \sum_{\chi \bmod q \atop \chi\neq\chi_0}
    c_\chi m_{\chi}(\rho),
\end{equation}
where $m_{\chi}(\rho)$ is the multiplicity of the zero $\rho$ for $L(s,\chi)$.

We believe that
\begin{equation}
    \label{eq:alpha_rho bound}
    \alpha_\rho = O_q(1)
\end{equation}
as $c_\chi=O(1)$, and we expect, for given $\chi$ and $\rho$,
$m_\chi(\rho)=0$ or $1$ (i.e. we expect $\rho$ to be at most a simple zero
of $L(s,\chi)$). More should be true-
distinct $L(s,\chi)$, for primitive $\chi$, presumably have distinct zeros.

For our purposes, we only need an estimate for the number of zeros in
an interval, and we will use the
following estimate:
\begin{equation}
    \label{eq:sum alpha_rho}
    \sum_{|\Im \rho|<T} |\alpha_\rho| = O_q(T \log{T}).
\end{equation}
The above bound follows from the asymptotic formula for the number
of zeros of $L(s,\chi)$ given in~\eqref{eq:N(T)}, which also implies
the following estimate that we will also use:
\begin{equation}
    \label{eq:sum alpha_rho 2}
    \sum_{n\leq|\Im \rho|<n+1} |\alpha_\rho| = O_q(\log{n}), \qquad \text{as $n \to \infty$}.
\end{equation}

\begin{lemma}
\label{lemma:alpha_rho nonzero}
Infinitely many $\alpha_\rho$ in~\eqref{eq:main formula zeros} are non-zero.
\end{lemma}
\begin{proof}
Notice that $A(x)$ has jump discontinuities at a relatively thin set of prime powers:
$R(x,1)$ and $R(x,q,a_j)$ jump when $x$ is a prime power $p^k$ with $k\geq 2$.
The remaining terms in $A(x)$ jump at prime powers $p^k$ with $p|q$, $k \geq
1$. Hence, overall, $A(x)$ has only finitely many jump discontinuities at the
primes, namely the primes $p$ that divide $q$. Furthermore every term that
appears in $B(x)$ is continuous with respect to $x$. But
$P(x)-\frac{r}{\phi(q)}\pi(x)$ is discontinuous at all primes (our assumption
that $r/\phi(q)<1$ enters here).

Therefore, the sum over $\rho$ in~\eqref{eq:main formula zeros}
must have infinitely many terms with $\alpha_\rho \neq 0$, otherwise the sum
over $\rho$ would be continuous a function for all $x$.

\end{proof}

\subsection{$\Omega$ results}

The fact that at least one $\alpha_\rho$ is non-zero
is a crucial point, and we are now in a position to obtain our
$\Omega$ results.

Let $\Theta$ be the lim sup of the real parts of the zeros $\rho$ such that
$\alpha_\rho \neq 0$, i.e. the zeros that appear in~\eqref{eq:main formula zeros}:
\begin{equation}
    \label{eq:Theta}
    \Theta = \lim\sup \{ \Re{\rho} | \alpha_\rho \neq 0\}.
\end{equation}
Equivalently, from~\eqref{eq:main formula},  $\Theta$ is the lim sup of the real parts of the poles of the function
\begin{equation}
    \label{eq:generating function a}
    -\sum_{\chi\neq \chi_0} c_\chi \frac{L'(s,\chi)}{L(s,\chi)},
\end{equation}
and also of the real parts of the singularities of
\begin{equation}
    \label{eq:generating function b}
    \sum_{\chi\neq \chi_0} c_\chi \log(L(s,\chi)).
\end{equation}
Notice that $\Theta \geq 1/2$, since the zeros of $L(s,\chi)$ that occur off the half line
(assuming GRH fails) come in pairs, $\rho$ and $1-\overline{\rho}$, symmetric about the line $\Re{s}=1/2$.

\subsection{$\Omega$ bound, assuming $\Theta>1/2$}
\label{sec:theta > 1/2}

We first assume that $\Theta>1/2$, i.e. that the GRH fails and that at least one zero
to the right of $\Re(s)=1/2$ survives in the explicit formula on taking the linear combination
in~\eqref{eq:main formula}.

We will prove the following Theorem.

\begin{theorem}
\label{thm:theta > 1/2}
Assume that $\Theta>1/2$. Then, for every $\delta>0$,
\begin{equation}
    \label{eq:lower bound a}
    P(x)- \frac{r}{\phi(q)} \pi(x) = \Omega (x^{\Theta-\delta}),
\end{equation}
with the implied constant depending on $\delta$ and $q$.
\end{theorem}

\begin{proof}
We do so by establishing, assuming $\Theta>1/2$, the estimate
\begin{equation}
    \label{eq:lower bound b}
    \sum_{j=1}^r \Pi(x,q,a_j) -\frac{r}{\phi(q)}\Pi(x) = \Omega (x^{\Theta-\delta}).
\end{equation}
The left hand side above is easier to work with than~\eqref{eq:lower bound a}, since
$L$-functions naturally count prime powers rather than just primes. Notice
that $\Pi(x,q,a)-\pi(x,q,a)=O_q(x^{1/2}/\log(x))$ and $\Pi(x) -
\pi(x)= O(x^{1/2}/\log(x))$ (see~\eqref{eq:R(x,1) estimate}
and~\eqref{eq:R(x,q,a) estimate} below). Hence, for
$0<\delta<\Theta-1/2$, we have that~\eqref{eq:lower bound b}
implies~\eqref{eq:lower bound a}. The bound~\eqref{eq:lower bound a} then holds
for every $\delta>0$, since taking $\delta$ larger gives a weaker bound.

Argue by contradiction. Assume that~\eqref{eq:lower bound b} does not hold,
i.e. that there exists a $\delta>0$ such that $|\sum_{j=1}^r \Pi(x,q,a_j) -\frac{r}{\phi(q)}\Pi(x)|
\ll x^{\Theta-\delta}$. Consider the Dirichlet integral (akin to Dirichlet
series, see Chapter 5 of~\cite{I}),
\begin{eqnarray}
    \label{eq:Dirichlet integral}
    &&\int_1^\infty \frac{\sum_{j=1}^r \Pi(x,q,a_j) -\frac{r}{\phi(q)}\Pi(x)}{x^{s+1}} dx \\
    &=& \frac{1}{s} \sum_{\chi} c_\chi \log(L(s,\chi)) - \frac{r}{s \phi(q)} \log(\zeta(s)) \notag \\
    &=& \frac{1}{s} \sum_{\chi\neq\chi_0} c_\chi \log(L(s,\chi)) + \frac{r}{s \phi(q)} \sum_{p|q} \log(1-p^{-s}). \notag
\end{eqnarray}
One can prove this identity, when $\Re{s}>1$, by observing that the numerator
of the integrand is a step function with steps at prime powers, and then integrating
termwise the contribution from each prime power. The assumption that $\Re{s}>1$ is used
to rearrange integration and summation and also to identify the resulting Dirichlet series
with the right hand side above.

Notice that the right hand side above has singularities (branch cuts)
coming from the zeros of $L(s,\chi)$, specifically from the $\rho$ with $\alpha_\rho \neq 0$,
i.e. those that survive the linear combination $\sum_{\chi\neq\chi_0} c_\chi \log(L(s,\chi))$.
There are also some additional singularities originating on the line $\Re{s}=0$.

Now, if the numerator of the above integrand  is $ \ll x^{\Theta-\delta}$, then
the left hand side of~\eqref{eq:Dirichlet integral} defines an analytic function for
$\Re{s}> \Theta-\delta$. But this contradicts, from the definition of $\Theta$,
the fact that the right hand side has singularities in this half plane. Therefore,
\begin{equation}
    \label{eq:omega result Pi}
    \sum_{j=1}^r \Pi(x,q,a_j) -\frac{r}{\phi(q)}\Pi(x) = \Omega(x^{\Theta-\delta}),
\end{equation}
for all $\delta>0$.

This establishes the Theorem. By taking $\delta$ sufficiently small, it also
proves Theorem~\ref{thm:1}, for residue classes, in the case that $\Theta>1/2$.

\end{proof}

\subsection{$\Omega$ estimate in the classical case, assuming $\Theta=1/2$}
\label{sec:theta = 1/2}

In this subsection, we assume
that $\Theta=1/2$. This can occur in two ways: either if the GRH holds, or if
the only zeros surviving the linear combination of explicit formulae arising
from~\eqref{eq:main formula} are on the half line.

\begin{theorem}
\label{thm:theta=1/2}
If $\Theta=1/2$, then
\begin{equation}
    \label{eq:omega bound on theta=1/2}
    P(x)- \frac{r}{\phi(q)} \pi(x) = \Omega(x^{1/2}/\log(x)).
\end{equation}
\end{theorem}

While we could modify the approach given in the previous subsection, it is complicated by the
presence of the squares of primes. We could adapt the approach described in
Ingham~\cite{I} for the problem of $\psi(x)-x$, and prove that
\begin{equation}
    \label{eq:Theta 1/2 and Ingham}
    \sum_{j=1}^r \Pi(x,q,a_j) -\frac{r}{\phi(q)}\Pi(x) = \Omega_{\pm}(x^{1/2}/\log{x}),
\end{equation}
i.e. that the difference of these prime {\it power} counting functions gets, in
size, as large as a constant times $x^{1/2}/\log{x}$, and points in {\em both}
positive and negative directions for infinite sequences of $x \to \infty$. Now
the squares of primes contribute an amount to $P(x)-\frac{r}{\phi(q)} \pi(x)$
that is asymptotically a constant times $x^{1/2}/\log(x)$, i.e. of the same
size as~\eqref{eq:Theta 1/2 and Ingham} , but always pointing in one direction.
Hence, estimate~\eqref{eq:Theta 1/2 and Ingham} would establish~\eqref{eq:omega
bound on theta=1/2}.

Instead, however, we will take an alternate approach that yields more
information. We will consider two mean square averages of the remainder term,
each giving a separate proof of~\eqref{eq:omega bound on theta=1/2}. Both
averages are of interest in their own right. To do so we first prove the following
Lemma which provides for a more manageable formula in comparison to Lemma~\ref{lemma:P explicit}.

\begin{lemma}
Let
\begin{equation}
    \label{eq:kappa}
    \kappa =
    \frac{1}{\phi(q)} \sum_{j=1}^{r}
    \sum_{b^2 = a_j \mod q} 1.
\end{equation}
Then, writing $\rho = 1/2 + i \gamma$, we have
\begin{equation}
    \label{eq:main estimate on GRH, finite sum nu}
    P(x) -\frac{r}{\phi(q)} \pi(x)
    = \frac{x^{1/2}}{\log{x}}
    \left( \sum_{0<|\gamma|< X} \alpha_\rho \frac{x^{i\gamma}}{\rho}
    + \nu +O\left(\frac{x^{1/2}\log(X)^2}{X}+\frac{1}{\log{x}}\right)\right),
\end{equation}
where $\nu$ is equal to $r/\phi(q)-\kappa$ plus, if the term $\rho=1/2$ appears
in~\eqref{eq:main formula zeros}, $2\alpha_{1/2}$.
\end{lemma}

\begin{proof}
We first bound each term that appears in $A(x)$, $B(x)$.
The prime number theorem and~\eqref{eq:R(x,1)} give
\begin{equation}
    \label{eq:R(x,1) estimate}
    R(x,1) =
    \pi(x^{1/2})/2 + O(x^{1/3}/\log{x})
    = x^{1/2}/\log{x} + O(x^{1/2}/\log(x)^2).
\end{equation}
Similarly, from the prime number theorem for arithmetic progressions, we have
\begin{equation}
    \label{eq:R(x,q,a) estimate}
    \sum_{j=1}^{r}
    R(x,q,a_j)
    = \kappa x^{1/2}/\log{x} + O_q(x^{1/2}/\log(x)^2),
\end{equation}
with the implied constant in the $O$ depending on $q$, and $\kappa$
defined above.

Finally, there are only finitely many $p|q$. Furthermore, $p^k \leq x$ implies
that $k \leq \log(x)/\log(p)$.
Hence $\sum_{p^k \leq x \atop  p \mid q} \log(p) = O_q(\log{x})$, and so
\begin{equation}
    \label{eq:A(x) finite p}
    -\frac{1}{2} \frac{\sum_{p^k \leq x \atop  p \mid q} \log(p)}{\log(x)}
    -\frac{1}{\log{x}} \sum_{ {\chi \mod q\atop \chi \neq \chi_0} \atop \chi \text{imprimitive}}
    c_\chi \sum_{p^k \leq x \atop  p | q} \log(p) \chi_1(p^k) =
    O_q(1).
\end{equation}
Putting these together gives
\begin{equation}
    \label{eq:A(x) estimate}
    A(x) = (r/\phi(q)-\kappa) x^{1/2}/\log{x} + O_q(x^{1/2}/\log(x)^2).
\end{equation}

To estimate $B(x)$, notice that
\begin{equation}
    \label{eq:B(x) estimate 1}
    \int_2^x \frac{ \sum_{p^k \leq t \atop  p \mid q} \log(p) }{t \log(t)^2} dt
    \ll_q
    \int_2^x \frac{dt}{t \log(t)} \ll \log\log{x}.
\end{equation}
Thus, because $\lambda$ and second line of~\eqref{eq:B(x)} are bounded, we have
\begin{eqnarray}
    \label{eq:B(x) estimate 2}
    B(x) \ll_q
    \left|
        \sum_{\chi \neq \chi_{0}} c_\chi
         \int_{2}^{x} \frac{\psi(t,\chi)}{t \log(t)^2} dt
    \right|
    + \log\log{x}.
\end{eqnarray}
Let
\begin{equation}
    G(x,\chi)= \int_{2}^{x} \psi(t,\chi) dt/t.
\end{equation}
Equation (5), page 117,
of Davenport~\cite{Davenport} gives the explicit formula with a rate of
convergence:
\begin{eqnarray}
    \label{eq:explicit formula b}
    \psi(x,\chi) &=& - \sum_{|\Im{\rho_\chi}|<X} \frac{x^{\rho_\chi}}{\rho_\chi}
    + O_q(x\log(xX)^2/X + \log{x}),
\end{eqnarray}
valid for $x \geq 2$ (or else for $x>1$ by adding a $1$ to the $O$-term).
Note that this formula, with the $O(\log{x})$ included, is true for both
primitive and imprimitive characters, and whether $x$ is equal to a prime power
or not. We have also absorbed the last three terms of~\eqref{eq:explicit
formula} into the $O$ term above. Thus, integrating and letting $X \rightarrow
\infty$, we have
\begin{equation}
    \label{eq:G}
    G(x,\chi)=- \sum_{\rho_{\chi}}
    \frac{x^{\rho_\chi}}{\rho_\chi^2}
    +O_q(\log(x)^2).
\end{equation}
The above series over $\rho_{\chi}$ converges absolutely
as can be seen from the asymptotic formula for
the number of zeros \cite[pg 101]{Davenport},
\begin{equation}
    \label{eq:N(T)}
    N(T,\chi) := \# \{\rho_{\chi}:  |\Im{\rho_{\chi}}| \leq T \} 
    = \frac{T}{\pi}\log\frac{qT}{2\pi} - \frac{T}{\pi} + O(\log T + \log q).
\end{equation}
Summing over $\chi$ we get
\begin{equation}
    \label{eq:G all chi}
    \sum_{\chi \neq \chi_{0}} c_\chi
    \int_{2}^{x} \frac{\psi(t,\chi)}{t} dt
    =- \sum_{\rho}
    \alpha_\rho \frac{x^{\rho}}{\rho^2}
    +O_q(\log(x)^2),
\end{equation}
with the sum over all non-trivial zeros $\rho$ of all $L(s,\chi)$ for
the modulus $q$.
The same coefficients $\alpha_\rho$ appear here as in~\eqref{eq:main formula zeros}
because the same linear combination of the terms involving the zeros $\rho$
appears as in the sum $\sum_{\chi \neq \chi_{0}} c_\chi \psi(x,\chi)$.

It follows, on integrating by parts, that
\begin{equation}
    \label{eq:bound on integral of psi}
    \sum_{\chi \neq \chi_{0}} c_\chi
    \int_{2}^{x} \frac{\psi(t,\chi)}{t \log(t)^2} dt
    \ll_q \frac{x^{1/2}}{\log(x)^2} \sum_{\rho}
    \frac{|\alpha_\rho|}{|\rho|^2} \ll_q \frac{x^{1/2}}{\log(x)^2}.
\end{equation}
The last bound follows, on summing by parts, from~\eqref{eq:sum alpha_rho}.

Thus, returning to~\eqref{eq:B(x) estimate 2}, we get
\begin{equation}
    \label{eq:B(x) estimate 3}
    B(x)
    \ll_q \frac{x^{1/2}}{\log(x)^2}.
\end{equation}

Thus, our estimates~\eqref{eq:A(x) estimate} and~\eqref{eq:B(x) estimate 3} for $A(x)$ and $B(x)$ give
\begin{equation}
    \label{eq:main estimate on GRH}
    P(x) - \frac{r}{\phi(q)} \pi(x)
    = \frac{1}{\log{x}}
    \left( \sum_{\rho} \alpha_\rho \frac{x^{\rho}}{\rho}
    + (r/\phi(q)-\kappa)x^{1/2} +O\left(x^{1/2}/\log{x}\right)\right),
\end{equation}
By~\eqref{eq:explicit formula b}, for $2 \leq x<X$,
we can write this as a finite sum over $\rho$, as expressed in
~\eqref{eq:main estimate on GRH, finite sum nu}.
The assumption $x<X$ is used here to simplify, in~\eqref{eq:explicit formula b},
$\log(xX)$ by $\log{X}$. We also use it below when estimating the contribution from the
above $O$ term.

Finally, we need to deal with the possibility of non-trivial zeros at $s=1/2$.
Such terms contribute $2 \alpha_{1/2} x^{1/2}/\log{x}$ to the sum
in~\eqref{eq:main estimate on GRH, finite sum nu}.
\end{proof}

\subsection{A mean square estimate of the average difference}
\label{sec:mean square a}

We continue with our proof of Theorem~\ref{thm:theta=1/2}.
Let
\begin{equation}
    \label{eq:Delta}
    \Delta(x):= \frac{\log{x}}{x^{1/2}} \left(P(x) - \frac{r}{\phi(q)} \pi(x)\right).
\end{equation}
Rather than work with $\Delta(x)$ directly, it is
technically easier to work with its average:
\begin{equation}
    \label{eq:average Delta}
    M(x) := \frac{1}{x} \int_2^x \Delta(t) dt =
    \sum_{\rho\neq 1/2} \alpha_\rho \frac{x^{i\gamma}}{\rho(i\gamma+1)}
    + \nu +O\left(1/\log{x}\right).
\end{equation}
The latter equality can be derived by integrating the bracketed expression
in~\eqref{eq:main estimate on GRH, finite sum nu} termwise, and letting $X \to
\infty$. Recall that $\nu$ is defined in~\eqref{eq:main estimate on GRH, finite
sum nu}

The Lemma below will be used to prove our $\Omega$ bound~\eqref{eq:omega bound on theta=1/2}.

\begin{lemma}
The following holds:
\begin{equation}
    \label{eq:mean square lemma}
    \lim_{Y \to \infty}
    \frac{1}{Y}\int_{\log{2}}^{Y} |M(e^y)|^2 dy =
    \sum_{0<|\gamma|} \frac{|\alpha_\rho|^2}{|\rho(i\gamma+1)|^2}
    +\nu^2 > 0.
\end{equation}

\end{lemma}

\begin{proof}
For any $\epsilon>0$, there exists $T=T(\epsilon)$ such that
\begin{equation}
    \label{eq:average Delta truncated}
    M(x) =
    \sum_{0<|\gamma|<T} \alpha_\rho \frac{x^{i\gamma}}{\rho(i\gamma+1)}
    + \nu + V(x),
\end{equation}
where
\begin{equation}
    \label{eq:V}
    V(x) < \epsilon,
\end{equation}
for all $x$ sufficiently large.
This can be obtained using estimate~\eqref{eq:sum alpha_rho 2} to show that the sum
in~\eqref{eq:average Delta} converges absolutely, and hence uniformly in $x$.

The natural scale at which to analyze the explicit formula is
logarithmic. Set $y=\log{x}$, and consider
\begin{equation}
    \label{eq:mean square a}
    \frac{1}{Y}\int_{\log{2}}^{Y} |M(e^y)|^2 dy.
\end{equation}
Substitute the right hand side of~\eqref{eq:average Delta truncated} for $M(e^y)$.
Now,
\begin{eqnarray}
    \label{eq:mean square main term}
    &&\lim_{Y \to \infty}
    \frac{1}{Y}\int_{\log{2}}^{Y}
    \left|
    \sum_{0<|\gamma|<T} \alpha_\rho \frac{e^{i\gamma y}}{\rho(i\gamma+1)}
        + \nu
    \right|^2 dy \notag \\
    &=&
    \sum_{0<|\gamma|<T} \frac{|\alpha_\rho|^2}{|\rho(i\gamma+1)|^2}
    +\nu^2,
\end{eqnarray}
which follows by multiplying
\begin{equation}
    \label{eq:sum rho}
    \sum_{0<|\gamma|<T} \alpha_\rho \frac{e^{i\gamma y}}{\rho(i\gamma+1)} + \nu
\end{equation}
by its conjugate, expanding, and noting that only the diagonal terms survive the limit
$Y\to\infty$. Next, the expression in~\eqref{eq:sum rho} is bounded for $y\in
\R$, and combining with~\eqref{eq:V} gives
\begin{eqnarray}
    \label{eq:mean square V}
    &&\frac{1}{Y}\int_{\log{2}}^{Y}
    \left(2\left|
    \sum_{0<|\gamma|<T} \alpha_\rho \frac{e^{i\gamma y}}{\rho(i\gamma+1)}
        + \nu
    \right| |V(e^y)| +
    |V(e^y)|^2\right) dy \notag \\
    &\ll& \epsilon+\epsilon^2,
\end{eqnarray}
for all $Y$ sufficiently large.

Since we may make $\epsilon$ as small as we wish, we get the equality expressed
in~\eqref{eq:mean square lemma}. We also have the inequality stated
in~\eqref{eq:mean square lemma} because, by Lemma~\ref{lemma:alpha_rho nonzero},
at least one $\alpha_\rho$ is non-zero.
\end{proof}

Hence,
\begin{equation}
    \label{eq:Omega M}
    M(e^y) = \Omega(1),
\end{equation}
i.e.
\begin{equation}
    \label{eq:Omega M b}
    M(x) = \Omega(1),
\end{equation}
which implies, from~\eqref{eq:average Delta}, that
\begin{equation}
    \label{eq:Omega Delta}
    \Delta(x) = \Omega(1),
\end{equation}
and hence from~\eqref{eq:Delta} we get~\eqref{eq:omega bound on theta=1/2}
of Theorem~\ref{thm:theta=1/2}.

\subsection{Unsmoothed mean square estimate}
\label{sec:mean square b}

In this subsection we give an alternate proof of the Omega
bound~\eqref{eq:omega bound on theta=1/2} by working out an unsmoothed mean
square.

Substituting $y=\log{x}$ and $Y=\log{X}$ in~\eqref{eq:main estimate on GRH,
finite sum nu}, we consider
\begin{equation}
    \label{eq:main estimate on GRH, finite sum a}
    P(e^y) - \frac{r}{\phi(q)} \pi(e^y)
    = \frac{e^{y/2}}{y} \left( \sum_{0<|\gamma|< e^Y} \alpha_\rho \frac{e^{i\gamma y }}{\rho}
    + \nu +O\left(\frac{e^{y/2} Y^2}{e^Y}+\frac{1}{y}\right)\right).
\end{equation}
Unlike $M(x)$ which was uniformly approximated by the finite
sum~\eqref{eq:average Delta truncated}, the above diverges absolutely and
cannot be uniformly approximated. However we can show
that we can approximate the sum with finitely many terms so that
the remainder term is uniformly small in mean square.

Thus, truncate the sum over $\rho$ at some large, but fixed $T$ (i.e.
independent of $y$), and consider, the mean square of the remainder. This is
essentially Lemma 2.2 of~\cite{RS}, but we provide slightly more detail here.
Thus, for $T \geq 1$, and $\log{2} \leq y$, we have
\begin{equation}
    \label{eq:main estimate on GRH, finite sum b}
    P(e^y) - \frac{r}{\phi(q)} \pi(e^y)
    = \frac{e^{y/2}}{y} \left( \sum_{0<|\gamma|< T} \alpha_\rho \frac{e^{i\gamma y }}{\rho}
    + \nu + r(y,T) \right),
\end{equation}
where, for all $Y$ satisfying $y \leq Y$,
\begin{equation}
    \label{eq:r}
    r(y,T)
    =  \sum_{T \leq |\gamma|<e^Y } \alpha_\rho \frac{e^{i\gamma y }}{\rho}
    + O\left(\frac{e^{y/2} Y^2}{e^Y}+\frac{1}{y}\right).
\end{equation}
The following Lemma gives a bound on the mean square of the remainder $r(y,T)$.
\begin{lemma}
\label{lemma:mean square}
Let  $T > 1$ and $Y>T^{1/2}/\log{T}$. Then,
\begin{equation}
    \frac{1}{Y/2} \int_{Y/2}^{Y} |r(y,T)|^{2} dy
    \ll_{q} \frac{\log(T)^2}{T}.
\end{equation}
\end{lemma}

\begin{proof}
Substitute \eqref{eq:r} into the integrand, and
use the inequality
\begin{equation}
    \label{eq:AGM}
    |a+b|^2 \leq 2(|a|^2 + |b|^2),
\end{equation}
which follows from the arithmetic geometric mean inequality $2|ab| \leq |a|^2 +|b|^2$,
to get
\begin{equation}
    \int_{Y/2}^{Y} | r(y,T) |^{2} dy
    \ll \int_{Y/2}^{Y} 
    \left|
        \sum_{T \leq |\gamma| \leq e^{Y}}
        \alpha_\rho \frac{e^{iy\gamma}}{1/2+i\gamma}
    \right|^{2} dy + 1/Y.
\end{equation}
The term $1/Y$ comes about from integrating the square of the $O$
term in~\eqref{eq:main estimate on GRH, finite sum a}. Multiplying out the
sum above by its conjugate, estimating the resulting integral, and extending the
double sum to infinity, the right hand side above becomes
\begin{eqnarray}
    &&
    \sum_{{T \leq |\gamma_1| \leq e^{Y}}\atop{T \leq |\gamma_2| \leq e^{Y}}}
    \frac{\alpha_{\rho_1} \bar{\alpha}_{\rho_2}}{\rho_1 \bar{\rho_2}}
    \int_{Y/2}^{Y} e^{iy(\gamma_1-\gamma_2)}
    dy + 1/Y \notag \\
    &\ll&
    \sum_{{T \leq |\gamma_1| \leq \infty}\atop{T \leq |\gamma_2| \leq \infty}}
    \frac{|\alpha_{\rho_1}||\alpha_{\rho_2}|}{|\rho_1| |\rho_2|}
    \min\left(Y,\frac{1}{|\gamma_1-\gamma_2|}\right) + 1/Y.
\end{eqnarray}

Breaking up the sum over zeros into unit intervals $|\gamma| \in [n,n+1)$, with $n \in \Z$, $n\geq T-1$,
and using~\eqref{eq:sum alpha_rho 2}, the above sum
is bounded by
\begin{equation}
    \label{eq:bound on double sum}
    \ll
    Y \sum_{n \geq T-1} \frac{\log(n)^2}{n^2}
    +\sum_{{n \geq T-1} \atop {m \geq n+1}} \frac{\log{m}}{m} \frac{\log{n}}{n}\frac{1}{m-n}.
\end{equation}
The first sum accounts for the contribution of the diagonal terms, i.e. where
$|\gamma_1|$ lies in an interval $[n,n+1)$ and $|\gamma_2|$ lies in $[m,m+1)$,
with $|m-n| \leq 1$. For such pairs of zeros, which could potentially be very
close, we use $Y$ as an upper bound for $\min(Y,1/|\gamma_1-\gamma_2|)$. For
all other pairs of zeros the quantity $1/|\gamma_1-\gamma_2| \ll 1/|m-n|$ is
much smaller than $Y$. This gives the second sum above, i.e. the off-diagonal
terms. We have also exploited symmetry in taking half the terms, i.e. we have
dropped $n \geq m+1$.

By comparing with the integral $\int_T^\infty \log(t)^2/t^2 dt$, we get,
on integrating by parts,
\begin{equation}
    \label{eq:bound diagonal}
    \sum_{n \geq T-1} \frac{\log(n)^2}{n^2} \ll \frac{\log(T)^2}{T}.
\end{equation}
To bound the off-diagonal contribution, break up the sum over $m$ into
the terms, $n+1 \leq m \leq 2n$, and the tail $m>2n$.
The first portion can be estimated as follows:
\begin{equation}
    \label{eq:bound off diagonal}
    \sum_{ {n \geq T-1} \atop {n+1 \leq m \leq 2n}} \frac{\log{m}}{m} \frac{\log{n}}{n}\frac{1}{m-n}
     \ll \sum_{n \geq T-1} \frac{\log(n)^2}{n^2} \sum_{n+1 \leq m \leq 2n} \frac{1}{m-n}
     \ll \sum_{n \geq T-1} \frac{\log(n)^3}{n^2} \ll \frac{\log(T)^3}{T}.
\end{equation}
For the contribution from the tail, use $1/|m-n| < 2/m$ when $m > 2n$:
\begin{equation}
    \label{eq:bound off diagonal b}
    \sum_{ {n \geq T-1} \atop {m > 2n}} \frac{\log{m}}{m} \frac{\log{n}}{n}\frac{1}{m-n}
     \ll
     \sum_{ n \geq T-1} \frac{\log{n}}{n} \sum_{m >2n} \frac{\log{m}}{m}\frac{1}{m}
     \ll \sum_{ n \geq T-1} \frac{\log(n)^2}{n^2} \ll \frac{\log(T)^2}{T},
\end{equation}
where we used $\sum_{m >2n} \log(m)/m^2 \ll \log(n)/n$ in passing from the second to third
expression.

Putting these bounds together gives
\begin{equation}
    \int_{Y/2}^{Y} |r(y,T)|^{2} dy
    \ll_{q} Y \frac{\log^{2} T}{T} +\frac{\log(T)^3}{T} + \frac{1}{Y}.
\end{equation}
For given $T$, and all $Y>T^{1/2}/\log{T}$, the first term on the right hand side dominates.
Dividing by $\frac{1}{Y/2}$ gives the Lemma.
\end{proof}

Returning to~\eqref{eq:main estimate on GRH, finite sum b}, we consider
the mean square:
\begin{equation}
    \label{eq:main mean square}
    \frac{1}{Y/2} \int_{Y/2}^Y 
    \left|(P(e^y) - \frac{r}{\phi(q)} \pi(e^y)) \frac{y}{e^{y/2}}\right|^2 dy
    =
    \frac{1}{Y/2} \int_{Y/2}^Y
        \left| \sum_{0<|\gamma|< T} \alpha_\rho \frac{e^{i\gamma y }}{\rho}
        + \nu + r(y,T) \right|^2 dy.
\end{equation}
The above equals
\begin{equation}
    \label{eq:main mean square b}
    \frac{1}{Y/2} \int_{Y/2}^Y
        \left| \sum_{0<|\gamma|< T} \alpha_\rho \frac{e^{i\gamma y }}{\rho}
        + \nu \right|^2 dy + E
\end{equation}
where
\begin{equation}
    \label{eq:E}
    |E| \ll
    \frac{1}{Y/2} \int_{Y/2}^Y
        \left| \sum_{0<|\gamma|< T} \alpha_\rho \frac{e^{i\gamma y }}{\rho}
        + \nu \right| |r(y,T)| dy
    + \frac{1}{Y/2} \int_{Y/2}^Y |r(y,T)|^2 dy.
\end{equation}
By multiplying the expression inside the absolute value of~\eqref{eq:main mean square b}
by its conjugate, and integrating termwise, we get
\begin{equation}
    \label{eq:main mean square c}
    \frac{1}{Y/2} \int_{Y/2}^Y
        \left| \sum_{0<|\gamma|< T} \alpha_\rho \frac{e^{i\gamma y }}{\rho}
        + \nu \right|^2 dy = \nu^2 + \sum_{0<|\gamma|< T} \frac{|\alpha_\rho|^2}{|\rho|^2}
        +O_T(1/Y).
\end{equation}
Next we estimate $E$. The bound~\eqref{eq:E} gives
\begin{equation}
    \label{eq:E 2}
    |E| \ll
    \max_{Y/2 \leq y \leq Y} \left| \sum_{0<|\gamma|< T} \alpha_\rho \frac{e^{i\gamma y }}{\rho} 
    + \nu \right|
    \frac{1}{Y/2} \int_{Y/2}^Y
        |r(y,T)| dy
    + \frac{1}{Y/2} \int_{Y/2}^Y |r(y,T)|^2 dy.
\end{equation}
Lemma~\ref{lemma:mean square} gives an estimate for the second integral:
\begin{equation}
     \frac{1}{Y/2} \int_{Y/2}^Y |r(y,T)|^2 dy \ll \frac{\log(T)^2}{T}.
\end{equation}
Now, from~\eqref{eq:sum alpha_rho 2}
\begin{equation}
    \left| \sum_{0<|\gamma|< T} \alpha_\rho \frac{e^{i\gamma y }}{\rho} + \nu \right|
    \ll \log(T)^2.
\end{equation}
Furthermore, the Cauchy-Schwarz inequality gives
\begin{equation}
    \frac{1}{Y/2} \int_{Y/2}^Y
        |r(y,T)| dy \ll
    + \frac{1}{Y/2} 
    \left(
        \int_{Y/2}^Y dy
        \int_{Y/2}^Y |r(y,T)|^2 dy
    \right)^{1/2},
\end{equation}
which, by Lemma~\ref{lemma:mean square}, is
\begin{equation}
    \ll \frac{\log{T}}{T^{1/2}}.
\end{equation}
Since we may choose $T$ as large as we please, we have
on combining the above estimates, that, as $Y \to \infty$,
\begin{equation}
    \label{eq:main mean square d}
    \frac{1}{Y/2} \int_{Y/2}^Y
    \left|(P(e^y) - \frac{r}{\phi(q)} \pi(e^y)) \frac{y}{e^{y/2}}\right|^2 dy
    \to
    \nu^2 + \sum_{\rho\neq 1/2} \frac{|\alpha_\rho|^2}{|\rho|^2}.
\end{equation}
As in the previous subsection, the right hand side above is positive because
at least one of the $\alpha_\rho$ is non-zero. Therefore,
\begin{equation}
    \left(P(e^y) - \frac{r}{\phi(q)} \pi(e^y)\right) \frac{y}{e^{y/2}} =\Omega(1),
\end{equation}
hence giving another proof of Theorem~\ref{thm:theta=1/2}.

Note that it is important that $r/\phi(q)<1$, since if we take $P$ to be the
set of all primes then it {\it can} be realized, in many ways, as a union
of primes in residue classes by taking all residue classes mod $q$, for any
positive integer $q$.
The reason the proof fails in this case is that $P(x)-\pi(x)$ is then
identically zero (and hence continuous),
giving a mean square for $P(x)-\pi(x)$, and more precisely,
of~\eqref{eq:main mean square d}, which is always zero.
The positivity of the right hand side of that equation requires there to be at least one
non-zero term appearing on the right hand side. However, from~\eqref{eq:c_chi},
all the $c_\chi$ and hence $\alpha_\rho$ are 0, and similarly for
the term $r/\phi(q)-\kappa$ which then equals $0$.

\setcounter{equation}{0}
\section{Generalization to Chebotarev sets}
\label{sec:chebotarev}

Here we generalize the problem to prime ideals and Chebotarev sets.

Therefore, let $L$ be a Galois extension of $K$ with Galois group
$G=\Gal(L/K)$. For a prime ideal $\p \in K$, we let the Artin symbol $(L/K,\p)$
denote the conjugacy class of Frobenius automorphisms corresponding to the
prime ideals $\Pri \in L$ that divide $\p$.

Given a conjugacy class $C$ of $G$, we let
\begin{equation}
    \label{eq:P C}
    \pi(L/K,C) =
    \left\{\p \in \pi(K) | \p
    \text{ unramified in $L$, $(L/K,\p) = C$}  \right\},
\end{equation}
consist of the unramified prime ideals $\p \in K$, and Frobenius conjugacy
class in $G$ equal to $C$. We define the counting function
\begin{equation}
    \label{eq:pi C}
    \pi(x,L/K,C) := \sum_{{N_{K/\Q}(\p) \leq x}\atop{\p \in \pi(L/K,C)}} 1
\end{equation}
to be the number of prime ideals in $\pi(L/K,C)$ with norm less than or equal to $x$.
Throughout what follows, we simply write $N\mathfrak{a}$ rather than $N_{K/\Q}(\mathfrak{a})$.

The Chebotarev density theorem states that
\begin{equation}
    \label{eq:chebotarev pnt}
    \pi(x,L/K,C) \sim \frac{|C|}{|G|} \Li(x),
\end{equation}
and the prime number theorem for prime ideals in $K$ states that
\begin{equation}
    \label{eq:pi(x,K)}
    \pi(x,K) := \{\p \in \pi(K)| N\p \leq x\} \sim \Li(x).
\end{equation}
Therefore, say we have a subset $P$ of prime ideals in $\pi(K)$ that is
realized, up to finitely many exceptions, as a finite union of Frobenius
conjugacy classes in the Galois group $G$ of some Galois extension $L$ of $K$.
We can restrict ourselves to the case of a single Galois extension $L$ for
similar reasons that we were able to restrict ourselves to a single modulus $q$
in the previous section. See the comments in the introduction in
Section~\ref{sec:Chebotarev intro}.

\subsection{Proof of Theorem 1 and 2}
\label{sec:proofs in general}

All the formulas used in the classical situation of residue classes in
Section~\ref{sec:classical} have analogues in the case of number fields. In
particular, the explicit formula for our situation has been worked out, with
remainder terms, by Lagarias and Odlyzko~\cite{LO}. We develop and collect
below the needed formulas.

Define
\begin{equation}
    \label{eq:psi C}
    \psi(x,L/K,C) := \sum_{ {{N\p^m \leq x}\atop{\text{$\p$ unramified}}}\atop{(L/K,\p)^m=C}}
                     \log N \p,
\end{equation}
\begin{eqnarray}
    \label{eq:Pi C}
    \Pi(x,L/K,C) &:=& \sum_{ {{N\p^m \leq x}\atop{\text{$\p$ unramified}}}\atop{(L/K,\p)^m=C}}
                      \frac{1}{m} \notag \\
                  &=& \pi(x,L/K,C) + R(x,L/K,C),
\end{eqnarray}
where 
\begin{eqnarray}
    \label{eq:R(x,L/K,C)}
    R(x,L/K,C) :=
               \sum_{ {{N\p^m \leq x}\atop{\text{$\p$ unramified,$m\geq 2$}}}\atop{(L/K,\p)^m=C}}
               \frac{1}{m},
\end{eqnarray}
so that
\begin{eqnarray}
    \label{eq:pi C vs Pi C}
    &&\pi(x,L/K,C) = \Pi(x,L/K,C) - R(x,L/K,C) \notag \\
    &=& \frac{\psi(x,L/K,C)}{\log(x)}
        + \int_2^x \frac{\psi(t,L/K,C)}{t \log(t)^2} dt -R(x,L/K,C).
\end{eqnarray}
Likewise, define
\begin{equation}
    \label{eq:Pi K}
    \Pi(x,K) := \sum_{ N\p^m \leq x} \frac{1}{m} = \pi(x,K) + R(x,K),
\end{equation}
where
\begin{eqnarray}
    \label{eq:R(x,K)}
    R(x,K) :=
               \sum_{ N\p^m \leq x \atop m\geq 2}
               \frac{1}{m}.
\end{eqnarray}
Thus,
\begin{eqnarray}
    \label{eq:pi K vs Pi K}
    \pi(x,K) &=& \Pi(x,K) - R(x,K) \notag \\
                 &=& \frac{\psi(x,K)}{\log(x)} 
                     + \int_2^x \frac{\psi(t,K)}{t \log(t)^2} dt -R(x,K).
\end{eqnarray}
We will also use
\begin{eqnarray}
    \label{eq:R(x,K) b}
    R(x,K) &=&
    \sum_{ N\p^2 \leq x}
    \frac{1}{2}
    +
    \sum_{ N\p^m \leq x \atop m\geq 3}
    \frac{1}{m} \notag \\
    &=& x^{1/2}/\log{x} + O(x^{1/3}/\log{x}),
\end{eqnarray}
which follows from the prime number theorem for ideals, with the implied
constant in the $O$ depending on $K$. Similarly, from the Chebotarev density
theorem, we have
\begin{equation}
    \label{eq:R C estimate}
    \sum_{j=1}^{r}
    R(x,L/K,C_j)
    = \kappa x^{1/2}/\log{x} + O(x^{1/2}/\log(x)^2),
\end{equation}
with the implied constant in the $O$ depending on $L/K$ and the $C_j$, and,
overriding the notation for $\kappa$ used earlier,
\begin{equation}
    \label{eq:kappa C}
    \kappa =
    \frac{1}{|G|} \sum_{j=1}^{r} |C_j|
    \sum_{b^2 \in C_j} 1,
\end{equation}
the inner sum counting the number of conjugacy class representatives $b \in G$ that,
when squared, lie in $C_j$.

To obtain an explicit formula for $\psi(x,L/K,C)$, Lagarias and Odlyzko mimic
the approach taken in Davenport for primes in arithmetic progression, using
the following linear combination of logarithmic derivatives of Artin $L$-functions
in order to extract primes ideals (and their powers) lying in the conjugacy class $C$:
\begin{eqnarray}
    \label{eq:F_C(s)}
    F_C(s) &:=& -\frac{|C|}{|G|} \sum_{\phi}\bar{\phi}(g) L'/L(s,\phi,L/K) \notag \\
    &=& \sum_{\p^m} \theta(\p^m) \log(N\p) (N\p)^{-ms},
\end{eqnarray}
where $g$ is any element of the conjugacy class $C$, $\phi$ runs over the irreducible
characters of $G$, and,
for unramified $\p$:
\begin{equation}
    \label{eq:theta is}
    \theta(\p^m) =
    \begin{cases}
        1 \quad \text{$(L/K,\p)^m = C$},\\
        0 \quad \text{otherwise},
    \end{cases}
\end{equation}
while, for ramified $\p$, $|\theta(\p^m)| \leq 1$. Notice that, while the right hand side
of~\eqref{eq:F_C(s)} resembles the Dirichlet series that gives the counting
function in~\eqref{eq:psi progressions}, there is a minor difference. Above,
and also in~\eqref{eq:F_C(s) b} below, the characters are primitive. The way to
interpret~\eqref{eq:psi progressions} so that it matches with the formula here,
is that each $\chi$ in~\eqref{eq:psi progressions} should be replaced by its
inducing character at a cost of $O(\log(x))$ to $\psi(x,q,a)$ coming from the
primes that ramify.

Brauer~\cite{B} proved that each Artin $L$-function can be written as a ratio
of Hecke $L$-functions, hence the linear combination of logarithmic derivatives
of Artin $L$-functions above can be written in terms of Hecke $L$-functions. In
our situation, the particular linear combination turns out, nicely, to have a
similar form to~\eqref{eq:F_C(s)}. Lagarias and Odlyzko use a construction
(Lemma 4.1 in their paper) of Deuring~\cite{D} to write
\begin{eqnarray}
    \label{eq:F_C(s) b}
    F_C(s) = -\frac{|C|}{|G|} \sum_{\chi}\bar{\chi}(g) L'/L(s,\chi,L/E),
\end{eqnarray}
where $\chi$ runs over the irreducible Hecke characters of $H=<g>$, the cyclic subgroup
generated by $g$, and $E$ is the fixed field of $H$.

The advantage of writing $F_C(s)$ in terms of Hecke characters is that the
analytic properties of Hecke $L$-functions are well established. Lagarias and
Odlyzko carry out a Perron integral in order to extract the Dirichlet
coefficients, with $N\p^m \leq x$, of $F_C(s)$.

Restricting to $2<x<X$, equation (7.4) of~\cite{LO} gives
\begin{equation}
    \label{eq:psi C estimate}
    \psi(x,L/K,C) = \frac{|C|}{|G|}
    \left(
        x - \sum_{\chi}\bar{\chi}(g)\sum_{|\Im \rho_\chi|<X} \frac{x^\rho_\chi}{\rho_\chi}
    \right) + \text{remainder}(x,X,L/K,C),
\end{equation}
where
\begin{equation}
    \text{remainder}(x,X,L/K,C)=
    O\left(\frac{x\log(X)^2}{X}+\log{x}\right),
\end{equation}
with the implied constant depending on $L/K$ and $C$. Here, $\rho$ runs over
all the non-trivial zeros of $L(s,\chi,L/E)$. The main term $\frac{|C|}{|G|} x$
arises from the principal character $\chi_0$ since $L(s,\chi_0,L/E)$ has, up to
finitely many Euler factors, $\zeta(s)$ as one of its factors, and hence a
simple pole at $s=1$.

Our remainder term is simpler than in (7.4) of Lagarias and Odlyzko
because we are taking $L/K$ to be fixed. Furthermore, $\text{remainder}(x,X,L/K,C)$
is a piecewise continuous function, with $O(\log{x})$ discontinuities at the points
$x=N\p^m$, where $\p$ runs over the ramified primes in $K$.

Substitute~\eqref{eq:psi C estimate} into~\eqref{eq:pi C vs Pi C}, apply the
estimate~\eqref{eq:R C estimate}, and then substitute all into
\begin{equation}
    \label{eq:P C from pi}
    P(x) = \lambda + \sum_{j=1}^r \pi(x,L/K,C_j),
\end{equation}
where $\lambda \in \Z$ accounts for finitely many exceptions, and
$x$ is sufficiently large (so that $x$ exceeds the norm of any of these
exceptions).
Letting
\begin{equation}
    \label{eq:sum C/G}
    \beta = \sum_{j=1}^r |C_j| / |G|,
\end{equation}
we have, on subtracting the analogous formula for $\beta \pi(x,K)$
and cancelling the main term coming from pole at $s=1$ of the factor 
of $\zeta(s)$ in the Dedekind zeta function $\zeta_K$, that, for $2 \leq x<X$,
and assuming GRH,
\begin{equation}
    \label{eq:main estimate C on GRH, finite sum}
    P(x) - \beta \pi(x,K)
    = \frac{x^{1/2}}{\log{x}} \left( \sum_{|\gamma|< X} \alpha_\rho \frac{x^{i\gamma}}{\rho}
    + (\beta-\kappa) +O\left(\frac{x^{1/2}\log(X)^2}{X}+\frac{1}{\log{x}}\right)\right),
\end{equation}
$\alpha_\rho \in \C$, and where the sum over $\rho$ is over the non-trivial
zeros of all relevant $L$-functions, namely the Hecke $L$-functions, for each
$C_j$ in~\eqref{eq:F_C(s) b}. More precisely,
\begin{equation}
    \label{eq:alpha_rho general}
    \alpha_{\rho} = - \frac{1}{|G|}
    \sum_{j=1}^r |C_j| \sum_{\chi \neq \chi_0} \bar{\chi}(g_j) m_{\chi}(\rho),
\end{equation}
with $m_{\chi}(\rho)$ the multiplicity of the zero $\rho$ for $L(s,\chi,L/E_j)$.

Bounds ~\eqref{eq:sum alpha_rho} and~\eqref{eq:sum alpha_rho 2} continue to hold, though
with the implied constants depending on $L$, $K$, and $C_j$ rather than on $q$.

As in the classical case, if the term $\rho=1/2$ appears in the sum, we absorb it into the
constant term. Thus, let $\mu=\beta-\kappa$ plus, in the event that $\rho=1/2$ appears
in the sum, $2\alpha_{1/2}$. The above then becomes
\begin{equation}
    \label{eq:main estimate C on GRH, finite sum b}
    P(x) - \beta \pi(x,K)
    = \frac{x^{1/2}}{\log{x}} \left( \sum_{0<|\gamma|< X} \alpha_\rho \frac{x^{i\gamma}}{\rho}
    + \mu +O\left(\frac{x^{1/2}\log(X)^2}{X}+\frac{1}{\log{x}}\right)\right),
\end{equation}

The cancellation of the main term deserves some
elaboration. The $L$-function corresponding to the principal character in
~\eqref{eq:F_C(s) b} factors as product of $\zeta_E(s)$ and Hecke
$L$-functions, and, because $K \subseteq E$, $\zeta_E(s)$ itself has
$\zeta_K(s)$ as a factor. The latter Dedekind zeta function is responsible for
cancellation in~\eqref{eq:main estimate C on GRH, finite sum} of the
main term. 
We therefore see that we could, in the statement of
Theorem 1 replace $\beta \pi(x,K)$ with
just $\beta \Li(x)$, since this
would have the same effect of cancelling the main term, with no further impact
on the form of the remaining terms. Similarly,
we could replace the counting function by $\frac{\beta}{\delta(Q)}Q(x)$ where
$Q \subseteq \pi(K)$ is another Chebotarev set with density $\delta(Q)$, or more generally,
take any difference of the form~\eqref{eq:F vs Li}, since the choice of $\eta$
there ensures cancellation of the main term.

Next, the jump discontinuities of the left hand side of~\eqref{eq:main estimate C on GRH,
finite sum b}, up to given $x$, outnumber those of the $O$ term of the right hand side, for
the same reason as in the case of residue classes mod $q$: the discontinuities
of the remainder term occur at $x=N\p^m$, $m\geq 2$, for $\p \in \pi(K)$,
coming from the terms $R(x,L/K,C_j)$ and $R(x,K)$, of which, because $m\geq 2$,
there are $O(x^{1/2}/\log{x})$ many. The other discontinuities come from the
ramified primes of which there are finitely many (and $O(\log{x})$ of their
powers, but these powers are already counted as discontinuities of $R(x,K)$).
On the other hand, the left hand side has jump discontinuities at all prime ideals $\pi
\in \pi(K)$ of which there are asymptotically $x/\log{x}$ many. We therefore
conclude, as previously, that infinitely many of the $\alpha_\rho$ must be
non-zero or else the sum over zeros would be a continuous function and the right hand side
would not have sufficiently many discontinuities.

Again this holds when we replace
$\beta \pi(x,K)$ with $\beta \Li(x)$,
or $Q(x)$ as above, though in the latter case we must also ensure that $Q$
does not essentially coincide with $P$, namely that the symmetric difference
$P \bigtriangleup Q$ is infinite. In each case the difference
between $P(x)$ and any of these counting functions has discontinuities
at a positive proportion of $N\p$ for primes ideals $\p \in \pi(K)$,
i.e. at $\gg x/\log{x}$ points. The same applies to the difference in~\eqref{eq:F vs Li}

As in the previous Section, define
\begin{equation}
    \label{eq:Theta K}
    \Theta = \lim\sup \{ \Re{\rho} | \alpha_\rho \neq 0\}.
\end{equation}
Then, $\Theta\geq 1/2$.

If $\Theta=1/2$, we have two mean square estimates analogous to those in
Sections~\ref{sec:mean square a} and~\ref{sec:mean square b}.
In order to carry out these estimates, we
also need a bound, as before, for the number of non-trivial zeros of a
Hecke $L$-function, $N(L,T) = |\{ \rho: L(\rho) = 0, |\Im{\rho}|\leq T, 0<
\Re{\rho}<1\}|$, in intervals of length one. Lagarias and Odlyzko prove, in Lemma 5.4 of~\cite{LO},
that $N(L,T+1)-N(L,T) = O(\log{T})$, with the implied constant depending on the
$L$-function, hence the method used to obtain the mean square estimate in the
case of Dirichlet $L$-functions and residue classes follows through as before
and we summarize the formulas.

Adapting the notation used in Section~\ref{sec:mean square a},
let
\begin{equation}
    \label{eq:Delta K}
    \Delta(x):= \frac{\log{x}}{x^{1/2}} \left(P(x) - \beta \pi(x,K)\right).
\end{equation}
and
\begin{equation}
    \label{eq:average Delta K}
    M(x) := \frac{1}{x} \int_2^x \Delta(t) dt.
\end{equation}
Then
\begin{eqnarray}
    \label{eq:mean square main term K}
    &&\lim_{Y \to \infty}
    \frac{1}{Y}\int_{\log{2}}^{Y}
    \left| M(e^y) \right|^2 dy \notag \\
    &=&
    \sum_{0<|\gamma|} \frac{|\alpha_\rho|^2}{|\rho(i\gamma+1)|^2}
    +\mu^2,
\end{eqnarray}
and
\begin{eqnarray}
    \label{eq:mean square b K}
    &&\lim_{Y\to \infty}
    \frac{1}{Y/2} \int_{Y/2}^Y
    \left|(P(e^y) - r \pi(e^y,K)) \frac{y}{e^{y/2}}\right|^2 dy \notag \\
    &=&
    \mu^2 + \sum_{\rho\neq 1/2} \frac{|\alpha_\rho|^2}{|\rho|^2}.
\end{eqnarray}
And, because at least one $\alpha_\rho$ is non-zero, both mean squares are positive.
From either, we can thus conclude as in Sections~\ref{sec:mean square a} or~\ref{sec:mean square b}
that
\begin{equation}
    P(x) - \beta \pi(x,K) = \Omega\left(\frac{x^{1/2}}{\log{x}} \right).
\end{equation}

This concludes our proof of Theorems~\ref{thm:1}-~\ref{thm:1c}.

We also get the following theorem, depending on the value of $\Theta$:
\begin{theorem}
\label{thm:3}
For every $\delta>0$,
\begin{equation}
    \label{eq:P(x) K Omega Theta>1/2}
    P(x) - \beta \pi(x,K) =
    \Omega(x^{\Theta-\delta}),
\end{equation}
with the implied constant in the $\Omega$ depending on $\delta$, and $P$.
\end{theorem}
\begin{proof}
If $\Theta=1/2$ then Theorem~\ref{thm:1} provides a stronger result and the
above therefore holds. If $\Theta>1/2$, the theorem follows, as in
Section~\ref{sec:theta > 1/2}, from the fact that at least one $\alpha_\rho$ is
non-zero and from the identity, initially derived with the assumption that
$\Re{s}>1$,
\begin{eqnarray}
    \label{eq:Dirichlet integral K}
    &&\int_1^\infty \frac{\sum_{j=1}^r \Pi(x,L/K,C_j) -\beta\Pi(x,K)}{x^{s+1}} dx \notag \\
    &=& \frac{1}{s} \sum_{j=1}^r \sum_{\chi_j} \bar{\chi}(g_j) \log(L(s,\chi_j,L/E_j))
    - \beta \log(\zeta_K(s)).
\end{eqnarray}
where, for $1\leq j\leq r$, $\chi_j$ runs over all the irreducible Hecke
characters of $H_j=<g_j>$, and $E_j$ is the fixed field of $H_j$.
\end{proof}

Finally, a similar theorem holds for a variety of
counting functions. As before, let $Q \subseteq \pi(K)$ be a Chebotarev set with the same
density $\beta$ as $P$, such that the symmetric difference $P \bigtriangleup Q$
is infinite, and $F$ any finite extension of $\Q$. Let $f(x)$ stand for
$\beta\Li(x)$, or $Q(x)$. For each such choice
of $f(x)$, the difference $P(x)-f(x)$ can be expressed as a linear combination of
explicit formulae having the same form as~\eqref{eq:main estimate C on GRH, finite sum}, though with the
constant term $\beta-\kappa$ replaced by $-\kappa$ in the case of
$f(x)=\beta\Li(x)$, and by $\kappa_2-\kappa$, when $f(x)=Q(x)$, where
$\kappa_2$ is the analogue of~\eqref{eq:kappa C} for the Chebotarev set $Q$.
Thus, defining $\Theta_f$ to be the analogue, for a given $f$, of~\eqref{eq:Theta K},
and likewise $\Theta_F$ for the difference in~\eqref{eq:F vs Li},
we have the following theorem:

\begin{theorem}
\label{thm:4}
Let $f(x)$ be as in the above paragraph, and $F(x)$ and $\eta$ be given
by~\eqref{eq:F} and \eqref{eq:lambda_0}. Then,
for every $\delta>0$,
\begin{equation}
    P(x) - f(x) =
    \Omega(x^{\Theta_f-\delta}),
\end{equation}
with the implied constant in the $\Omega$ depending on $\delta$, $P$, and $f$, and
\begin{equation}
    F(x) - \eta \Li(x) =
    \Omega(x^{\Theta_F-\delta}),
\end{equation}
with the implied constants in the $\Omega$ depending on $\Theta$, and $\eta$. The latter result
also holds (with $\Theta_F$ adjusted accordingly), if we replace
$\Li(x)$ by $\pi(x,K)$, with the restriction that not all $\eta_j$ are equal.
\end{theorem}

{\bf Acknowledgements:} We would like to thank Barry Mazur for raising the problem addressed
in this paper. Jeffrey Lagarias and Jean-Pierre Serre provided some feedback. The referee
provided helpful comments.

\end{document}